\documentclass[12pt]{amsart}

\usepackage{amsmath,amssymb,amscd,amsxtra}
\usepackage{latexsym}
\usepackage[dvips]{graphics,epsfig}

\newcommand{\am}[3]{W\,(L^{#1}, \ell_{#2}^{#3}\,)}

\newcommand{\CHI}{\hbox{\raise .3ex \hbox{$\chi$}}}

\newcommand{\ip}[2]{\langle#1,#2\rangle}
\newcommand{\bigip}[2]{\bigl\langle #1, \; #2 \bigl\rangle}

\newcommand{\nm}[1]{\|#1\|}

\newcommand{\R}{{\mathbb R}}
\newcommand{\Sc}{{\mathcal{S}}}

\newcommand{\Z}{{\mathbb Z}}
\newcommand{\M}[1]{\mathcal{M}^{#1}}

\newtheorem{theorem}{Theorem}[section]
\newtheorem{cor}[theorem]{Corollary}
\newtheorem{lemma}[theorem]{Lemma}
\newtheorem{prop}[theorem]{Proposition}
\theoremstyle{definition}
\newtheorem{definition}[theorem]{Definition}
\newtheorem{example}[theorem]{Example}
\newtheorem{remark}[theorem]{Remark}
\newtheorem{conjecture}[theorem]{Conjecture}

\newcommand{\begeq}{\begin {equation}}
\newcommand{\eq}{\end{equation}}
\newcommand{\bs}{\begin {split}}
\newcommand{\es}{\end{split}}
\newcommand{\bp}{\begin {prop}}
\newcommand{\ep}{\end {prop}}
\newcommand{\bt}{\begin {theorem}}
\newcommand{\et}{\end {theorem}}
\newcommand{\bc}{\begin {cor}}
\newcommand{\ec}{\end {cor}}
\newcommand{\bl}{\begin {lemma}}
\newcommand{\el}{\end {lemma}}
\newcommand{\bpf}{\begin {proof}}
\newcommand{\epf}{\end {proof}}
\newcommand{\bi}{\begin {itemize}}
\newcommand{\ei}{\end {itemize}}
\newcommand{\ben}{\begin {enumerate}}
\newcommand{\een}{\end {enumerate}}
\newcommand{\brem}{\begin {remark}}
\newcommand{\erem}{\end {remark}}
\renewcommand{\kappa}{\wp}

\newcommand{\norm}[1]{\left\|{#1} \right\|}

\newcommand{\la}{\langle}
\newcommand{\ra}{\rangle}

\newcommand{\cB}{\mathcal B}

\newcommand{\HH}{{\mathcal H}}

\newcommand{\PP}{{I\kern-.3em P}}
\newcommand{\ZZ}{{\mathbb Z}}
\newcommand{\TT}{{\mathbb T}}
\newcommand{\RR}{{I\kern-.3em R}}
\newcommand{\CC}{{I\kern-.6em C}}
\newcommand{\NN}{{I\kern-.3em N}}
\newcommand{\GG}{\ZZ^d}
\newcommand{\GGG}{\mathbb G}

\begin{document}

\title{\bf Invertibility of the Gabor frame operator \\ on the Wiener amalgam space}

\author {Ilya~A.~Krishtal}

\address{Department of Mathematical Sciences \\
Northern Illinois University \\
DeKalb, IL 60115, USA}

\email{krishtal@math.niu.edu}

\author{Kasso~A.~Okoudjou}

\address{Department of Mathematics \\
University of Maryland\\
College Park, MD 20742, USA}

\email{kasso@math.umd.edu}


\date{\today}


\begin{abstract} We use a generalization of Wiener's $1/f$ theorem to prove that for a Gabor frame with the generator in the Wiener amalgam space $\am{\infty}{}{1}(\R^{d})$, the corresponding frame operator is invertible on this space. Therefore, for such a Gabor frame, the generator of the canonical dual belongs also to $\am{\infty}{}{1}(\R^{d})$. 
\end{abstract}

\maketitle \pagestyle{myheadings} \thispagestyle{plain}
\markboth{I. A. KRISHTAL AND K. A. OKOUDJOU}{INVERTIBILTY OF THE GABOR FRAME OPERATOR}

\section{Introduction}\label{intro}
For $\alpha , \beta>0$ and $g \in L^{2}(\R^{d})$, let  $M_{\beta m}g(x) = e^{-2\pi i\beta m\cdot x}g(x)$ and $T_{\alpha n}g(x) = g(x-\alpha n)$. The collection $\mathcal{G}(g, \alpha, \beta) = \{g_{m,n} = M_{\beta m}T_{\alpha n}g,\ m,n\in\Z^d\}\subset L^{2}(\R^{d})$
is called a Gabor frame if there exist positive constants $0<A\leq B< \infty$ such that for each $f\in L^{2}(\R^{d})$ the following inequalities hold:
\begin{equation}\label{frameineq}
A\, \nm{f}_{L^{2}}^{2} \leq \sum_{m, n \in \Z^{d}} |\ip{f}{g_{m,n}}|^{2} \leq B\, \nm{f}_{L^{2}}^{2}.
\end{equation}
Equivalently, the frame condition can be restated in terms of the \emph{frame operator} 
$S_g:  L^{2}(\R^{d})\to L^{2}(\R^{d})$
associated to $\mathcal{G}(g, \alpha, \beta)$ and given by  
\begeq \label{gfo}
S_g f = \sum_{m,n\in\Z^d} \la f, g_{m,n}\ra g_{m,n},\quad f\in L^{2}(\R^{d}). 
\eq

In particular, $\mathcal{G}(g, \alpha, \beta)$ is a frame for $L^2$ if and only if $$A \nm{f}_{L^{2}}^{2} \leq \ip{S_{g}f}{f} \leq B \nm{f}_{L^{2}}^{2},
\quad\mbox{for all }  f\in L^{2}(\R^{d}).$$ 
When $\mathcal{G}(g, \alpha, \beta)$ is a frame for $L^2$ then with $\tilde g = S_g^{-1} g$ the following reconstruction formulas hold 
$$f=\sum_{m,n\in\Z^{d}} \ip{f}{\tilde{g}_{m,n}}g_{m,n} = \sum_{m,n \in \Z^{d}}\ip{f}{g_{m,n}}\tilde{g}_{m,n}.$$ Moreover,  $\mathcal{G}(\tilde{g}, \alpha, \beta)$ is also a Gabor frame for $L^2$ called the \emph{canonical dual} (Gabor) frame. We refer to \cite{Dau92, Gr01, hw89} for more on Gabor analysis. 

A central question in Gabor analysis is to find conditions on $\alpha, \beta>0$ and $g \in L^2$ such that $\mathcal{G}(g, \alpha, \beta)$ is a Gabor frame for $L^2$. Moreover, in many applications it is desirable to find a Gabor frame such that the generator $g$ and its canonical dual $\tilde{g}$ have the same properties, e.g., same type of decay and/or smoothness. For example, it was proved that if $g \in \Sc(\R^{d})$ the Schwartz class, and $\alpha, \beta >0$ are such that $\mathcal{G}(g, \alpha, \beta)$ is a Gabor frame, then $\tilde{g} \in \Sc$ \cite{jan95}. Similarly, let $\phi(x)=e^{-x^2}$ and define the short time Fourier transform of a tempered distribution $f \in \Sc'$  by $$V_{\phi}f(x, \omega)= \ip{f}{M_{\omega}T_{x}\phi}=\int_{\R^{d}}f(t)\, \overline{\phi(t-x)}\, e^{-2\pi i x\cdot \omega}\, dt.$$ The space $\M{1}(\R^{d})$ of all $f\in \Sc'$ such that $$\nm{f}_{\M{1}}=\iint_{\R^{2d}}|V_{\phi}f(x, \omega)|\, dx\, d\omega <\infty$$ is known as the Feichtinger algebra \cite{Fei81}.  In this context, Gr\"ochenig and Leinart proved a deep result that shows that if  $ g\in \M{1}(\R^{d})$ and $\alpha, \beta >0$ are such that $\mathcal{G}(g, \alpha, \beta)$ is a Gabor frame, then $\tilde{g} \in \M{1}$ \cite{grle03}. More specifically, using the so-called Janssen's representation of the Gabor frame operator, which converges absolutely in $\cB(L^{2}(\R^{d}))$ whenever $g \in \M{1}$, Gr\"ochenig and Leinart, recasted the question into a non-commutative version of the celebrated $1/f$ Wiener's lemma (\cite{Wie32}) involving the twisted convolution. We refer to \cite{fol89} for background on the twisted convolution.

Note that similarly to $\Sc$,  $\M{1}$ is also invariant under the Fourier transform \cite{Fei81, Gr01}. Moreover, it is trivially seen that  $\mathcal{G}(g, \alpha, \beta)$ is a Gabor frame if and only if $\mathcal{G}(\hat{g}, \beta, \alpha)$ is a Gabor frame. Therefore, the above results simply say that the generator of a Gabor frame and its canonical dual have the same time-frequency concentration. Furthermore, results involving time or frequency only conditions were proved in \cite{bo, bojan}. Finally, we refer to \cite{CCJ, Cho} for some related results.

In this paper, we prove a similar result that shows that if a Gabor frame $\mathcal{G}(g, \alpha, \beta)$ is generated by $g \in \am{\infty}{\nu}{1}$ then the generator of its canonical dual $\tilde{g} \in \am{\infty}{\nu}{1}$, where $\nu$ is an admissible weight (see Definition \ref{wegt}) and the space $\am{\infty}{\nu}{1}$ is a weighted Wiener amalgam space. Some particular cases of  our results were first obtained by Walnut with some extra  conditions \cite{wal90}. The unweighted amalgam space $\am{\infty}{}{1}$ was introduced by N.~Wiener in connection with the Tauberian Theorems \cite{Wie32}. The precise definition of $\am{\infty}{\nu}{1}$ is given as follows.  

For $\alpha>0$ let $Q_{\alpha}=[0, \alpha)^d$ and $\chi_{Q_{\alpha}}$ be the characteristic function of $Q_{\alpha}$. Let also $\nu:\GG\to[1,\infty)$
be a weight function. A function
$f \in \am{\infty}{\nu}{1}$ if and only if 
\begin{equation}\label{defwiama}
\nm{f}_{\am{\infty}{\nu}{1}}= \sum_{n\in \Z^{d}}\nm{f\cdot T_{n\alpha}\chi_{Q_{\alpha}}}_{L^{\infty}}\nu(n) < \infty.
\end{equation}
Moreover, equipped with this norm $\am{\infty}{\nu}{1}$ is a Banach space, whose definition is independent of $\alpha$ in the sense of equivalent norm. Furthermore, the following embeddings can be easily established: $\Sc \subset \M{1} \subset \am{\infty}{\nu}{1}\subset \am{\infty}{}{1}\subset L^2.$

We wish to point out that, the condition $g \in \am{\infty}{}{1}$ is not enough for the Gabor frame operator to admit an absolutely convergent Janssen's representation (\cite[Section 7.2]{Gr01}). Therefore, our result does not follow from \cite{grle03}. To prove our result, we rely instead on another representation of the Gabor frame operator: the Walnut's representation \cite{wal90}. In particular, our proof is derived using this representation of the frame operator together with a far-reaching generalization of the $1/f$ Wiener's lemma due to Baskakov \cite{Bas97', Bas97}. 
This particular extension of the 
lemma turns out to be the most suitable for us among it's numerous analogs, see \cite{Bal, Bas90, GKW, grle03, Jaf, Ku90, Sjo, Qiyu, Shu}, etc. Another analog that is suitable for us and pertains to the almost periodic situation will appear in \cite{BalKri}.
We observe that Kurbatov \cite{Ku99a, Ku99b} seems to be  the first to use this type of result in the context of amalgam spaces. We also refer to \cite{fastro} for relevant results. 

Our paper is organized as follows. In Section~\ref{pre} we state the precise version of the $1/f$ Wiener's lemma that is suited to our result. Moreover, we introduce the two main tools used in proving our result: the Walnut's representation of the Gabor frame operator, and the bracket product. In Section~\ref{mai} we state and prove our main result and, furthermore, outline a second and different proof. The second approach, however, relies on a conjecture that we have not yet been able to prove.

\section{Preliminaries}\label{pre}

\subsection{Wiener's lemma for Fourier series of operators.}

In this section we present a reformulation of a non-commutative Wiener's lemma proved by Baskakov in \cite{Bas97', Bas97}. 
We begin 
by introducing a notion of a Fourier series of a linear operator with respect to a representation of a compact Abelian group.

Although the results of this section hold for an arbitrary compact Abelian group $\GGG$, we restrict our attention to 
$\GGG = \TT_{1/\beta}^{d} \simeq Q_{1/\beta} =[0,1/\beta)^d$,
which is specifically tailored for our application. We use an additive form for the group operation on $\GGG$.
For an isometric strongly continuous representation $U: \GGG \to \cB(X)$, where $\cB(X)$ is the algebra of all bounded linear operators
on a (complex) Banach space $X$, we define 
$\tilde U: \GGG \to \cB(\cB(X))$ via
\[\tilde U(\theta) A = U(\theta)AU(-\theta),\quad \theta\in\GGG,\ A\in\cB(X).\]
Following \cite{Bas97'}, the Fourier series of an operator $A\in\cB(X)$ with respect to the representation $U$ is, by definition, the 
Fourier series of the function $\hat A: \GGG\to \cB(X)$ given by
\[\hat A(\theta) = \tilde U(\theta)A, \quad \theta\in\GGG.\]
Recall that this Fourier series is 
\begin{equation}
\hat A(\theta)f \simeq \sum_{k\in\GG} e^{2\pi i \beta \theta\cdot k} A_k f, \quad f\in X,
\end{equation}
and the Fourier coefficients $A_k \in \cB(X)$ are given by
\begin{equation}
A_k f = \beta^{d} 
\int_{Q_{1/\beta}} e^{-2\pi i \beta \theta\cdot k}\hat A(\theta)f d\theta
= \int_{[0,1)^d} e^{-2\pi i \theta\cdot k}\hat A(\frac\theta\beta)f d\theta. 
\end{equation}
Observe that the Fourier coefficients are eigen-vectors of the corresponding representation, i.e.
\begeq\label{ev}
\tilde U(\theta) A_k = e^{2\pi i \beta\theta\cdot k} A_k,\ k\in\GG.
\eq

\begin{example}\label{spec}
We are especially interested in the case when $X = \HH$ is the 
Hilbert space $L^2(\R^d)$ 
and the representation $U:\GGG\to\cB(\HH)$ is defined by
\[U(\theta)f(x) = M_{\beta\theta}f(x) = e^{2\pi i \beta\theta\cdot x}f(x),\quad f\in\HH.\]
Then the operators $T_{\frac n\beta}$, $n\in\GG$, are eigen-vectors for the representation $\tilde U$:
\[\tilde{U}T_{\tfrac{n}{\beta}}f(x)=U(\theta)T_{\frac n\beta}U(-\theta)f(x) = e^{2\pi i \beta \theta\cdot n} f(x - \frac n\beta),\quad f\in\HH,\ n\in \ZZ^d.\] 
This implies that any operator $A\in \cB(\HH)$ has Fourier coefficients of the form
\[A_n = \mathfrak{G}_n T_{\frac n\beta}, \quad n\in\GG,\] 
where $\mathfrak{G}_n$ commute with $M_{\beta\theta}$ for all $\theta \in Q_{1/\beta}$. Therefore, and following \cite{Ku99a},
$\mathfrak{G}_n$ is an operator of multiplication by a uniquely determined function $G_n\in L^\infty(\R^d)$ and 
\[\norm{A_n} = \norm{\mathfrak{G}_n} = \norm{G_n}_{L^\infty},\quad n\in\GG.\]

Observe that this construction remains partly valid when $\HH$ is a closed subspace of the Hilbert space $L^2(\R^d)$ invariant with respect to modulations
$M_\theta$, $\theta\in\R$, and translations $T_{\frac n\beta}$, $n\in\GG$. In this case, however, we can no longer infer that the operators $\mathfrak{G}_n$, $n\in\GG$, \emph{uniquely} determine a function $G_n\in L^\infty(\R^d)$ such that $(\mathfrak{G}_nf)(x) = G_n(x)f(x)$. Hence, we can only guarantee that
\[\norm{A_n} = \norm{\mathfrak{G}_n}\le \norm{G_n}_{L^\infty},\quad n\in\GG.\]
\end{example}

As usual when Wiener's lemma is discussed, we are interested in linear operators whose Fourier series
are summable or summable with a weight.

\begin{definition} \label{wegt}
An \emph{admissible weight} is a function $\nu: \GG\to [1,\infty)$ such that
\begin{enumerate}
\item $\nu$ is an even function, that is, $\nu(-n) = \nu(n)$, for all $n\in\GG$,
\item $\nu(k+n)\leq \nu(k)\nu(n)$, for all $k,n\in\GG$, and
\item $\lim\limits_{k\to\infty}k^{-1}\ln\nu(kn)=0$, for all $n\in\GG$.
\end{enumerate}
\end{definition}

For an admissible weight $\nu$, we consider a Banach algebra
\[\cB_{\nu}(X) = \{A\in \cB(X):\ \sum_{k\in\GG} \norm{A_k}\nu(k) <\infty\}\]
of linear operators with $\nu$-summable Fourier series. If $\nu \equiv 1$ we get the algebra
$\cB_{1}$ of operators with summable Fourier series.

The result below follows immediately from \cite[Theorem 2 and Remark]{Bas97'} or \cite{BalKri}.

\begin{theorem}\label{wie}
Let $\nu$ be an admissible weight and $A\in \cB_{\nu}(X)$ be an invertible operator.
Then $A^{-1}\in \cB_{\nu}(X)$. In particular, if $A$ is invertible and $A\in\cB_1$, then  $A^{-1}\in\cB_1$.
\end{theorem}

The next corollary, which follows immediately from Theorem \ref{wie} and Example \ref{spec},
plays a key role in establishing our main result.

\begin{cor}\label{abkwiener}
Let $\beta >0$, $\nu$ be an admissible weight, and $\HH=L^2(\R^d)$.  Assume that an invertible operator
$S \in \mathcal{B}(\HH)$
has a $\nu$-summable Fourier series, that is,
$$
S  = \sum_{n\in\Z^{d}} \mathfrak{G}_{n} T_{\frac n\beta}=\sum_{n \in \Z^{d}}G_{n} \cdot T_{\frac n\beta},
$$  
where $\mathfrak{G}_{n} \in \mathcal{B}(\HH)$ is the operator of multiplication by a function 
$G_n\in L^\infty(\R^d)$ and  $\sum\limits_{n\in\Z^{d}} \nm{\mathfrak{G}_{n} T_{\frac n\beta}} \nu(n) = \sum\limits_{n\in \Z^{d}} \nm{{G}_{n}}_{L^\infty} \nu(n)< \infty$. 
Then the inverse operator $S^{-1} \in \mathcal{B}(\HH)$ has also a $\nu$-summable Fourier series, that is, there exists a sequence of functions 
$\tilde{G}_{n} \in L^{\infty}(\R^d)$ such that $$S^{-1}=\sum_{n \in \Z^{d}}\mathfrak{\tilde{G}}_{n}T_{\frac n\beta}=\sum_{n \in \Z^{d}}\tilde{G}_{n}\cdot T_{\frac n\beta}, \mbox{ and }\sum_{n \in \Z^{d}} \nm{\tilde{G}_{n}}_{L^\infty} \nu(n)< \infty.$$
\end{cor}

The following result is a different version of Corollary~\ref{abkwiener} that deals with operators defined on closed subspaces of $L^{2}$. 

\begin{cor}\label{abkwiener2}
Let $\alpha >0$, $\nu$ be an admissible weight, and $\HH$ be a closed subspace of the Hilbert space $L^2(\R^d)$ invariant with respect to modulations
$M_\theta$, $\theta\in\R$, and translations $T_{\alpha n}$, $n\in\GG$. Assume that an invertible operator
$S \in \mathcal{B}(\HH)$ 
has a $\nu$-summable Fourier series, that is,
$$
S  = \sum_{n \in \Z^{d}} \mathfrak{G}_{n} T_{\alpha n} \quad \mbox{and}\quad \sum\limits_{n\in\Z^{d}} \nm{\mathfrak{G}_{n}}\nu(n) < \infty, 
$$  
where the operators $\mathfrak{G}_{n} \in \mathcal{B}(\HH)$ 
commute with $M_{\frac\theta\alpha}$ for all $\theta \in Q_{\alpha}$.
Then the inverse operator $S^{-1} \in \mathcal{B}(\HH)$ also has a $\nu$-summable Fourier series, that is, there exists a sequence of 
operators $\tilde{\mathfrak{G}}_{n} \in \mathcal{B}(\HH)$ commuting with $M_{\frac\theta\alpha}$ for all $\theta \in Q_{\alpha}$
and such that 
 $$S^{-1}=\sum_{n \in \Z^{d}}\mathfrak{\tilde{G}}_{n}T_{\alpha n}\quad \mbox{and}\quad
 \sum_{n\in\Z^{d}} \nm{\tilde{\mathfrak{G}}_{n}}\nu(n) < \infty.$$
\end{cor}

\begin{remark}\label{problem}
Observe that in Corollary \ref{abkwiener2}, even though  
$\mathfrak{G}_{n},\tilde{\mathfrak{G}}_{n} \in \mathcal{B}(\HH)$ are operators of multiplication by functions 
$G_n,\tilde G_n\in L^\infty(\R^d)$, respectively,
 we can no longer guarantee that
\[\sum_{n\in \Z^{d}} \nm{G_{n}}_{L^{\infty}} \nu(n)< \infty \mbox{ and }
\sum_{n \in \Z^{d}} \nm{\tilde{G}_{n}}_{L^{\infty}} \nu(n)< \infty.\]
\end{remark}

In the following theorem we use Corollary \ref{abkwiener} to establish the boundedness of operators $S\in\cB_\nu(\HH)$ on the Wiener amalgam space $\am{\infty}{\nu}{1}$. 

\begin{theorem}\label{forbound}
Let $\beta >0$, $\nu$ be an admissible weight, and $\HH=L^2(\R^d)$.  Assume that
$S= \sum_{n\in\Z^{d}} \mathfrak{G}_{n} T_{\frac n\beta}\in\cB_{\nu}(\HH)$ has a $\nu$-summable Fourier series.
Then $S$ defines a bounded operator from $\am{\infty}{\nu}{1}$ to $\am{\infty}{\nu}{1}$.
\end{theorem}

\begin{proof}
Since $S\in\cB_\nu$, 
$$
S  = \sum_{n \in \Z^{d}} \mathfrak{G}_{n} T_{\frac n\beta}=\sum_{n \in \Z^{d}}G_{n} \cdot T_{\frac n\beta},
$$  
where $\mathfrak{G}_{n} \in \mathcal{B}(\HH)$, $n\in\Z^d$, are the operators of multiplication by  functions 
$G_n\in L^\infty(\R^d)$ and  $\sum\limits_{n\in\Z^{d}} \nm{\mathfrak{G}_{n} T_{\frac n\beta}} \nu(n) = \sum\limits_{n\in \Z^{d}} \nm{{G}_{n}}_{L^\infty} \nu(n)< \infty$. 
Hence,
\begin{align*}
\nm{S f}_{\am{\infty}{\nu}{1}} & =\norm{\sum_{n\in\GG}(T_{n/\beta}f)G_{n}}_{\am{\infty}{\nu}{1}}\\
& = \sum_{k\in\GG}\norm{\sum_{n\in\GG}\chi_{k+Q_{1/\beta}}(T_{n/\beta}f)G_{n}}_\infty \nu(k)\\
& \le \sum_{k,n\in\GG}\norm{\chi_{k+Q_{1/\beta}}(T_{n/\beta}f)}_\infty \nu(k) \norm{G_{n}}_\infty \\
& \le \sum_{k,n\in\GG}\norm{\chi_{n+k+Q_{1/\beta}}f}_\infty \nu(n+k)\norm{G_{n}}_\infty \nu(n)\\
& \le  \norm{f}_{\am{\infty}{\nu}{1}}\sum_{n\in\GG} \norm{G_{n}}_\infty \nu(n) \le \infty,
\end{align*}
and the proof is complete.
\end{proof}

\subsection{Bracket product and Walnut representation.} 
Unless stated otherwise, in all that follows we assume that $\alpha, \beta >0$, $\nu$ is an admissible weight, $g \in \am{\infty}{\nu}{1}$ is such that $\mathcal{G}(g, \alpha, \beta)$ is a Gabor frame for $L^2$, and $\tilde{g}=S_{g}^{-1}g \in L^2$ is the generator of the canonical dual frame. 

Let us recall a few properties of the bracket product widely used in the study of shift invariant systems \cite{bdvr2, cala, rs97}.

For $f,h \in L^2(\R^d)$ and $\alpha > 0$ the $\alpha$-\emph{bracket product} of $f$ and $h$ is the $\alpha-$periodic function,
 which is a periodization of $f\cdot \bar{h} \in L^{1}(\R^{d})$:
\begin{equation}\label{brackpro}
[f, h]_{\alpha}(x) = \sum_{k \in \Z^{d}} (f\cdot \overline{h})(x - \alpha k),
\end{equation} 
$x \in Q_{\alpha}=[0, \alpha)^{d}$. Note that the series in~\eqref{brackpro} converges for $a. e. \, x \in Q_{\alpha}$. Observe also that the (formal) Fourier series of this periodic function is
\begeq
[f,h]_\alpha(x)\simeq \alpha^{-d}\sum_{n\in\GG} \la f, M_{\tfrac{n}{\alpha}} h\ra e^{\tfrac{2\pi i n\cdot x}{\alpha} }.
\eq

For $g \in \am{\infty}{\nu}{1}$, let
\begin{equation}\label{gperiod}
G_{n}(x)= [g, T_{\tfrac{n}{\beta}}g]_{\alpha}(x) = \sum_{k\in \Z^{d}}(g\cdot T_{\tfrac{n}{\beta}}\overline{g})(x - \alpha k). 
\end{equation} 

For the dual generator $\tilde{ g}=S_{g}^{-1}g$ we let $\tilde{G}_{n}$ be
\begin{equation}\label{tildegn}
\tilde{G}_{n} (x)=[\tilde g, T_{\tfrac{n}{\beta}}\tilde g]_{\alpha}(x) = \sum_{k \in \Z^{d}} (\tilde{g} \cdot T_{\tfrac{n}{\beta}} \overline{\tilde{g}})(x - \alpha k),
\end{equation}
which is well-defined for $a. e.\,   x \in Q_{\alpha}$ since $\tilde{g}\cdot T_{\tfrac{n}{\beta}}\overline{\tilde{g}} \in L^{1}$.

It can be shown  \cite[Lemma 6.3.1]{Gr01}, see also \cite[Lemma 5.2]{GHO} and \cite[Lemma 2.1]{wal90}, that if $g \in \am{\infty}{\nu}{1}$, there exists a constant $C$ which depends only on $\alpha$, $\beta$  and $d$ such that

\begin{equation}\label{linfty1}
\sum_{n\in \Z^{d}}\nm{G_{n}}_{L^{\infty}(Q_{\alpha})} \nu(n)\leq C \nm{g}_{\am{\infty}{\nu}{1}}^{2} < \infty.
\end{equation}
and
\begin{equation}\label{linfty1bis}
\sum_{n\in \Z^{d}}\norm{ [g, T_{\alpha n}g]_{1/\beta}}_{L^{\infty}(Q_{1/\beta })} \nu(n)\leq C \nm{g}_{\am{\infty}{\nu}{1}}^{2} < \infty.
\end{equation}

In Lemma \ref{bound} below, we shall present a converse to the above statement that will play a key role in obtaining our main result. To prove it we introduce
the \emph{Walnut representation} of the Gabor frame operators $S_{g}$ and $S_{\tilde g}$.
Following \cite[Proposition 7.1.1]{Gr01}, we have
\begin{equation}\label{weakwalnut1}
\bigip{S_{g} f}{h} =\bigip{\beta^{-d}\sum_{n \in \Z^{d}} {G}_{n} \cdot T_{\tfrac{n}{\beta}}f}{h}, 
\end{equation} 
\begin{equation}\label{weakwalnut}
\bigip{S_{\tilde{g}} f}{h} =\bigip{\beta^{-d}\sum_{n \in \Z^{d}} \tilde{G}_{n} \cdot T_{\tfrac{n}{\beta}}f}{h},
\end{equation} 
for all $f, h$ bounded and compactly supported. In fact, because $g\in \am{\infty}{\nu}{1}$ the operator $S_g$ has a \emph{strong} Walnut representation \cite[Theorem 6.3.2]{Gr01}
\begin{equation}\label{walnutrep}
S_{g}f = \beta^{-d}\sum_{n\in\Z^d} \mathfrak{G}_{n} T_{\tfrac{n}{\beta}}f=\beta^{-d}\sum_{n \in \Z^{d}}G_{n} \cdot T_{\tfrac{n}{\beta}}f,\quad f\in L^2(\R^d),
\end{equation}  
where $\mathfrak{G}_{n} \in \mathcal{B}(L^{2}(\R^{d}))$ is the operator corresponding to multiplication by the bounded function $G_n$ given in~\eqref{gperiod}. Moreover, observe that $\nm{\mathfrak{G}_{n}}=\nm{G_{n}}_{L^{\infty}(Q_{\alpha})}$ and, therefore, \eqref{linfty1} implies $\sum_{n}\nm{\mathfrak{G}_{n}}_{op} \nu(n)= \sum_{n}\nm{G_{n}}_{L^{\infty}(Q_{\alpha})} \nu(n)< \infty$.

By Example \ref{spec} and the above inequality, the Walnut representation \eqref{walnutrep} is the Fourier series of $S_g\in\cB_\nu(L^2(\R^d))$ 
with respect to the representation
\[U(\theta)f(x) = e^{2\pi i \beta \theta\cdot x}f(x),\quad \theta\in Q_{1/\beta},\ f\in L^2(\R^d).\]
Similarly, the weak Walnut representation \eqref{weakwalnut} implies that
\[S_{\tilde g} \simeq \beta^{-d}\sum_{n\in\GG}\widetilde {\mathfrak G}_n T_{\tfrac{n}{\beta}},\]
where $\widetilde {\mathfrak G}_n$ is the operator of multiplication by $\tilde G_n$, $n\in\GG$,
is the Fourier series of $S_{\tilde g}$.

\begin{lemma}\label{bound}
Assume that $g\in L^2(\R^d)$ is such that $\mathcal{G}(g, \alpha, \beta)$ is a Gabor frame for $L^2$ and
\begin{equation}\label{condi}
\sum_{n\in\GG}\nm{G_{n}}_{L^{\infty}(Q_{\alpha})} \nu(n)< \infty.
\end{equation}
Then the frame operator $S_g$ is a bounded operator from $\am{\infty}{\nu}{1}$ to $\am{\infty}{\nu}{1}$.
\end{lemma}

\begin{proof}
Since \eqref{condi} implies $S_g\in\cB_\nu$, the result follows immediately from Theorem \ref{forbound}.
\end{proof}

\begin{remark}
For a different proof of Lemma~\ref{bound} we refer to \cite{GHO}, and \cite[Theorem 3.1]{wal90}.
A similar result is also proved in \cite[Theorem 7.2]{CCJ}. Following that proof it can be shown that condition~\eqref{condi} in Lemma \ref{bound} is not only sufficient but also necessary.
\end{remark}

\section{Main results}\label{mai}
We are now ready to prove our main result.

\begin{theorem}\label{main1}
Let $\alpha, \beta >0$, $g \in \am{\infty}{\nu}{1}$, be such that $\mathcal{G}(g, \alpha, \beta)$ and $\mathcal{G}(\tilde g, \alpha, \beta)$ be canonical dual Gabor frames for $L^{2}(\R^{d})$, where $\tilde{g}=S_{g}^{-1}g$.
Then $\tilde g \in \am{\infty}{\nu}{1}$.
\end{theorem}

\begin{proof}
As mentioned above, since $g \in \am{\infty}{\nu}{1}$, the Gabor frame operator $S_g$ has a $\nu$-summable Fourier series \eqref{walnutrep}.
By Corollary \ref{abkwiener},
its inverse $S_{g}^{-1}=S_{\tilde{g}}$, \cite[Lemma 5.1.6]{Gr01},  also has a $\nu$-summable Fourier series, i.e.,
\[ S_{\tilde g}= \beta^{-d}\sum_{n\in \Z^{d}} \widetilde{\mathfrak{G}}_{n} T_{\tfrac{n}{\beta}},\] 
and
\[\sum_{n\in\GG}\nm{\tilde{G}_{n}}_{L^{\infty}(Q_{\alpha})} \nu(n)< \infty.\]
It remains to apply Lemma \ref{bound} to conclude that $\tilde{g}\in \am{\infty}{\nu}{1}$.
\end{proof}

\begin{remark} It is known that for a Gabor frame $\mathcal{G}(g, \alpha, \beta)$ of $L^{2}(\R^{d})$ the system $\mathcal{G}(\gamma, \alpha, \beta)$ is a dual frame if and only if $\gamma = \tilde{g} +h$ where $\tilde{g}$ is the canonical dual, and $h \in L^{2}(\R^{d})$ is such that $\ip{h}{M_{\tfrac{n}{\alpha}}T_{\tfrac{m}{\beta}}g}=0$ for all $m, n \in \Z^{d}$, e.g., see, \cite[Lemma 7.6.1]{Gr01}, and \cite{wr}. Therefore, in view of Theorem~\ref{main1}, one can ask whether for a given Gabor frame $\mathcal{G}(g, \alpha, \beta)$ with $g \in \am{\infty}{}{1}$, all the dual frames belong to the same space. In general, the answer to this question is no, as shown by the following example.

Let $g = \chi_{[0,1]}$ be the generator of the Gabor frame $\mathcal{G}(g, 1/2, 1)$ and $\tilde g$ be the generator of its canonical dual. 
For an arbitrary sequence $(a_k)\in\ell^2\backslash \ell^1$, define $h\in L^2$ via
\[h(x) = \sum_{k\in\Z} a_k \chi_{[k, k+1)}(x)e^{2\pi i x}.\] 
Then $\la h, M_{2m}T_n g\ra = 0$ for all $m,n\in\Z$ and, therefore, $\tilde{g}+h$ is a dual generator for $g$.
However, by construction, $h \notin \am{\infty}{}{1}$, and hence $\tilde{g}+h \notin \am{\infty}{}{1}$. 
\end{remark}

We recall that if $\mathcal{G}(g, \alpha, \beta)$ is a Gabor frame for $L^{2}(\R^{d})$, then $\mathcal{G}(g^{\dagger}, \alpha, \beta)$ is a tight frame, where
 $g^{\dagger}=S_{g}^{-1/2}g.$ Indeed, for all $f\in L^{2}(\R^{d})$, we have $$f=S_{g}^{-1/2}S_{g}S_{g}^{-1/2}f = \sum_{k, l}\ip{f}{g_{k,l}^{\dagger}} g_{k,l}^{\dagger}.$$ The next result proves that if $g\in \am{\infty}{\nu}{1}$, then $g^{\dagger}\in \am{\infty}{\nu}{1}$. More specifically we have

\begin{cor}\label{sqroot}
Let $\alpha, \beta >0$ and $g \in \am{\infty}{\nu}{1}$ be such that $\mathcal{G}(g, \alpha, \beta)$ is a Gabor frame for $L^{2}(\R^{d})$. Then $g^{\dagger} \in \am{\infty}{\nu}{1}$.
\end{cor}

\begin{proof}
Let $g \in \am{\infty}{\nu}{1}$ and $\alpha, \beta >0$ be such that $\mathcal{G}(g, \alpha, \beta)$ is a Gabor frame for $L^{2}(\R^{d})$.  Note that $S_{g} \in \cB_{\nu}(L^{2}(\R^{d}))$ is a positive definite operator on $L^{2}(\R^{d})$. Therefore, we can use the Riesz-Dunford functional calculus \cite[Ch. VII]{DSch} to get
\[S_g^{-1/2} = \frac1{2\pi i}\int_\Gamma \lambda^{-1/2}(S_g-\lambda I)^{-1}d\lambda,\]
where $\Gamma$ is a positively oriented contour in the right complex half-plane surrounding the spectrum of $S_g$. By Corollary \ref{abkwiener} 
the above integral converges in the norm of $\cB_\nu$ and we get
$S_{g}^{-1/2} \in \cB_{\nu}(L^{2}(\R^{d}))$. Consequently, Theorem~\ref{forbound} can be used to conclude that 
 $g^{\dagger}=S_{g}^{-1/2}g \in \am{\infty}{\nu}{1}.$
\end{proof}

\begin{remark} As mentioned in the Introduction, Theorem~\ref{main1} and Corollary~\ref{sqroot} were first proved under extra assumptions on $g$, $\alpha$ and $\beta$ in \cite[Corollary 3.5]{wal90}.
\end{remark}

We wish to conclude by outlining a {\emph {possible}} alternative approach to the proof of Theorem~\ref{main1} which does not use Lemma \ref{bound}. Instead, this proof relies upon the following propositions, the first of which, to our knowledge has not been proved before despite its simplicity.

\begin{prop}\label{keylemma1}
Let $g \in \am{\infty}{\nu}{1}$ and 
$\mathcal{G}(g, \alpha, \beta)$ and $\mathcal{G}(\tilde g, \alpha, \beta)$ be canonical dual Gabor frames for $L^{2}(\R^{d})$.
Then for almost every $x \in Q_{\tfrac{1}{\beta}}$ and all $k \in \Z^d$
\begin{equation}\label{convo}
[\tilde {g}, T_{\alpha k}g]_{\tfrac{1}{\beta}} (x) = \beta^{-d}\sum_{n \in \Z^{d}} \overline{[g, T_{\alpha n}g]_{\tfrac{1}{\beta}}}(x-\alpha k)\, [\tilde {g}, T_{\alpha(k +n)}\tilde {g}]_{\tfrac{1}{\beta}}(x)
\end{equation}
and we have the following norm estimates
\begin{equation}\label{convest}
\sum_{k \in \Z^{d}} \nm{[\tilde{g}, T_{\alpha k}g]_{\tfrac{1}{\beta}}}_{L^{\infty}}\nu(k) \leq \beta^{-d}\, \sum_{n\in \Z^{d}} \nm{[g, T_{\alpha n}g]_{\tfrac{1}{\beta}}}_{L^{\infty}} \nu(n)\, \sum_{n \in \Z^{d}} \nm{[\tilde{g}, T_{\alpha n}\tilde{g}]_{\tfrac{1}{\beta}}}_{L^{\infty}}\nu(n),
\end{equation} whenever the right-hand side is finite.
\end{prop}

\begin{remark}
From Theorem \ref{main1} and \eqref{linfty1bis}, we know that the right-hand side of \eqref{convest} is always finite under assumptions of
Proposition \ref{keylemma1}. However, we cannot use this fact if we want to give an alternative proof of Theorem \ref{main1}.
\end{remark}

\begin{proof}
Because $\mathcal{G}(g, \alpha, \beta)$ and $\mathcal{G}(\tilde{g}, \alpha, \beta)$ are dual frames we have $$g = \sum_{r, s \in \Z^{d}}\bigip{g}{M_{\beta s}T_{\alpha r}g} M_{\beta s}T_{\alpha r}\tilde{g}$$ and so $$ T_{\alpha k}g = \sum_{r, s \in \Z^{d}}\bigip{g}{M_{\beta s}T_{\alpha r}g}\, e^{-2\pi i \alpha \beta s\cdot k}\, M_{\beta s}T_{\alpha (r +k)}\tilde{g}.$$ Therefore, for $a. e.\, x \in Q_{\tfrac{1}{\beta}}$

\begin{align*}
[\tilde{g}, T_{\alpha k}g]_{\tfrac{1}{\beta}}(x) & = \sum_{n \in \Z^{d}} (\tilde{g}\cdot T_{\alpha k}\overline{g})(x - \tfrac{n}{\beta})\\
& = \sum_{r, s \in \Z^{d}} \overline{\bigip{g}{M_{\beta s}T_{\alpha r}g}}\, e^{2 \pi i \alpha \beta s \cdot k}\, e^{-2\pi i \beta s \cdot x}\, [\tilde{g}, T_{\alpha(k+r)}\tilde{g}]_{\tfrac{1}{\beta}}(x)\\
& = \sum_{r \in \Z^{d}} [\tilde{g}, T_{\alpha(k+r)}\tilde{g}]_{\tfrac{1}{\beta}}(x)\sum_{s\in \Z^{d}} \overline{\bigip{g}{M_{\beta s}T_{\alpha r}g}}\, e^{2 \pi i  \beta s \cdot (x -\alpha k)}\\
& = \beta^{-d}\, \sum_{r \in \Z^{d}} \overline{[g, T_{\alpha r}g]_{\tfrac{1}{\beta}}}(x-\alpha k)\, [\tilde{g}, T_{\alpha(k+r)}\tilde{g}]_{\tfrac{1}{\beta}}(x),
\end{align*}
where the last equation follows from Carleson's theorem  since for $g\in \am{\infty}{}{1}$ the Fourier series of 
$[g, T_{\alpha r}g]_{\tfrac{1}{\beta}}(x-\alpha k) = \beta^d\, \sum_{s\in \Z^{d}} \overline{\bigip{g}{M_{\beta s}T_{\alpha r}g}}\, e^{2 \pi i  \beta s \cdot (x -\alpha k)}$ with $L^2$ convergence.
Consequently, $$ \nm{[\tilde{g}, T_{\alpha k}g]_{\tfrac{1}{\beta}}}_{L^{\infty}} \leq \beta^{-d}\, \sum_{r \in \Z^{d}} \nm{[g, T_{\alpha r}g]_{\tfrac{1}{\beta}}}_{L^{\infty}}\, \nm{[\tilde{g}, T_{\alpha (k+r)}\tilde{g}]_{\tfrac{1}{\beta}}}_{L^{\infty}},$$ $$\sum_{k \in \Z^{d}} \nm{[\tilde{g}, T_{\alpha k}g]_{\tfrac{1}{\beta}}}_{L^{\infty}} \nu(k)\leq \beta^{-d}\, \sum_{r \in \Z^{d}} \nm{[g, T_{\alpha r}g]_{\tfrac{1}{\beta}}}_{L^{\infty}}\nu(r)\, \sum_{r \in \Z^{d}}\nm{[\tilde{g}, T_{\alpha r}\tilde{g}]_{\tfrac{1}{\beta}}}_{L^{\infty}}\nu(r),$$
and the proof is complete. 
\end{proof}

\begin{prop}\label{ab}
Let $g \in \am{\infty}{\nu}{1}$ and $\alpha, \beta >0$ be such that $\mathcal{G}(g, \alpha, \beta)$ is a Gabor frame for $L^{2}(\R^{d})$. Assume that $\mathcal{G}(\tilde g, \alpha, \beta)$ is its canonical dual  frame. Then  
$\sum_{r \in \Z^{d}}\nm{[\tilde{g}, T_{\frac r\beta}\tilde{g}]_{\alpha}}_{L^{\infty}}\nu(r) < \infty.$ 
\end{prop}

\begin{proof}
The result follows immediately from Corollary~\ref{abkwiener}, just as in the proof of Theorem \ref{main1}.
\end{proof}

Note that Proposition~\ref{ab} implies that $$\norm{[\tilde{g}, \tilde{g}]_{\alpha}}_{L^{\infty}} = \norm{\sum_{k\in \Z^{d}} |\tilde{g}(\cdot- \alpha k)|^{2}}_{L^{\infty}} \leq B < \infty $$ and it follows that $\tilde{g} \in L^{\infty}(\R^{d}).$

\begin{conjecture}\label{conj1}
Let $\alpha,\beta>0$ and $\tilde g\in L^2(\R^d)$ be such that $$\sum_{r \in \Z^{d}}\norm { [\tilde{g}, T_{\tfrac{r}{\beta}} \tilde{g}]_{\alpha}}_{L^{\infty}}\nu(r) < \infty.$$
We conjecture that in this case
\[\sum_{n \in \Z^{d}}\norm{[\tilde{g}, T_{\alpha n} \tilde{g}]_{\tfrac{1}{\beta}}}_{L^{\infty}} \nu(n)< \infty.\]
\end{conjecture}

Remark \ref{problem} indicates a major obstacle in proving
Conjecture~\ref{conj1}. However, if the conjecture is true, we can give the proof of Theorem~\ref{main1} as follows. 

Let $m_{k}(x) = \beta^{-d}[\tilde{g}, T_{\alpha k}g]_{\tfrac{1}{\beta}}(x)$, $k\in\Z^d$. These functions are well-defined because $g, \tilde{g} \in L^{2}(\R^{d})$ and so $ \tilde{g}\cdot T_{\alpha k}\overline{g} \in L^{1}(\R^{d})$.  Consequently, $m_{k} \in L^{\infty}(Q_{\tfrac{1}{\beta}})$ and is $\tfrac{1}{\beta}-$periodic.  
Moreover, $\hat{m}_{k}(l) = C_{g}\tilde{g}(\alpha k, \beta l)$, $k, l\in \Z^d$, are the Gabor coefficients of $\tilde{g}$ with respect to the frame $\mathcal{G}(g, \alpha, \beta)$ \cite{GHO}. By Proposition~\ref{keylemma1} and Conjecture \ref{conj1},
$$\sum_{k\in \Z^{d}}\nm{m_{k}}_{L^{\infty}} \nu(k)\leq \beta^{-d}\, \sum_{n\in \Z^{d}} \nm{[g, T_{\alpha n}g]_{\tfrac{1}{\beta}}}_{L^{\infty}} \nu(n)\, \sum_{n \in \Z^{d}} \nm{[\tilde{g}, T_{\alpha n}\tilde{g}]_{\tfrac{1}{\beta}}}_{L^{\infty}}\nu(n)< \infty.$$ Hence, it follows from  \cite[Theorem 4.4]{GHO} that 
$\tilde{g}\in \am{\infty}{\nu}{1}$.

\section{Acknowledgements}

We would like to thank R.~Balan, A.~Baskakov, P.~Casazza, O.~Christensen, B.~Farrell, K.~Gr\"ochenig, C.~Heil, A.J.E.M.~Janssen, and T.~Strohmer for their useful comments. We are also grateful to J.J.~Benedetto for organizing the February Fourier Talks at the University of Maryland, where we started this project.

\end{document}